\documentclass[A4, 12pt]{article}
\usepackage[top=1in,bottom=1in,left=1in,right=1in]{geometry}

\usepackage[sumlimits,intlimits]{amsmath}
\usepackage{amsthm,amssymb,amsfonts}
\usepackage{graphicx,color}
\usepackage{bm}
\usepackage{color}

\usepackage[
linktocpage=true,
colorlinks=true,
linkcolor=blue,
citecolor=blue,
urlcolor=blue
]{hyperref}
\usepackage[numbers,sort&compress]{natbib}

\usepackage{epstopdf}

\allowdisplaybreaks

\newtheorem{definition}{Definition}[section]
\newtheorem{proposition}[definition]{Proposition}
\newtheorem{lemma}[definition]{Lemma}
\newtheorem{theorem}[definition]{Theorem}

\newtheorem*{remark*}{Remark}

\theoremstyle{definition}

\title{Conditional Stability for Single Interior Measurement}
\author{
Naofumi \textsc{Honda},
\footnote{Department of Mathematics, Hokkaido University, Faculty of
Science, 060-0810, Japan.
\newline
$\qquad$E-mail: honda@math.sci.hokudai.ac.jp.
\newline
$\qquad$Supported in part by JSPS Grant-in-Aid for Scientific Research (C) (23540178) of the Japan Society for the Promotion of Science.
}
\and
Joyce \textsc{McLaughlin},
\footnote{Department of Mathematical Sciences
Rensselaer Polytechnic Institute, 12180, USA.
\newline
$\qquad$E-mail: mclauj@rpi.edu.
\newline
$\qquad$Supported by Office of Naval Research grants NOOO14-13-1-0388 and NOOO14-08-0432.
}
\and
Gen \textsc{Nakamura}
\footnote{Department of Mathematics, Inha University, 402-751, South Korea.
\newline
$\qquad$E-mail: nakamuragenn@gmail.com.
\newline
$\qquad$Supported in part by JSPS grants-in-aid Grant-in-Aid for Scientific Research (B) (22340023) of the Japan Society for the Promotion of Science when he was affiliated in Department of Mathematics, Hokkaido University, 060-0810, Japan.
}
}

\begin{document}

\maketitle

\begin{abstract}
An inverse problem to identify unknown coefficients of a partial differential equation by a single interior measurement is considered. The equation considered in this paper is a strongly elliptic second order scalar equation which can have complex coefficients in a bounded domain with $C^2$ boundary and single interior measurement means that we know a given solution of the equation in this domain. The equation includes some model equations arising from acoustics, viscoelasticity and hydrology. We assume that the coefficients are piecewise analytic. Our major result is the local H\"older stability estimate for identifying the unknown coefficients. If the unknown coefficients is a complex coefficient in the principal part of the equation, we assumed a condition which we named admissibility assumption for the real part and imaginary part of the difference of the two complex coefficients. This admissibility assumption is automatically satisfied if the complex coefficients are real valued. For identifying either the real coefficient in the principal part or the coefficient of the 0-th order of the equation, the major result implies the global uniqueness for the identification.
\end{abstract}



\pagestyle{myheadings}
\thispagestyle{plain}
\markboth{
NAOFUMI HONDA, JOYCE MCLAUGHLIN AND GEN NAKAMURA
}{CONDITIONAL STABILITY FOR SINGLE INTERIOR MEASUREMENT}

\newcommand{\Anabla}{\nabla\hspace{-3pt}{}_{A}\,}
\newcommand{\CAnabla}{\nabla}
\newcommand{\Unabla}{\Anabla u_1 \cdot \overline{\nabla u_1}}

\section{Introduction}
\noindent In order to make our description of the background more concise we first introduce two assumptions, and formulate the forward and inverse
problems.

\medskip
\noindent {\it Assumption 1:}
 Let $\Omega\subset\mathbb{R}^n$ $(2\le n\in\mathbb{N})$ be
a bounded domain with a $C^2$ smooth boundary $\partial\Omega$. Also, let
$A(x)$ be a strictly positive Hermitian matrix on $\overline\Omega$ with entries in $C^1(\overline\Omega)$
and $\gamma\in C^1(\overline\Omega)$ be either a real valued positive or complex valued function with positive real and imaginary parts on $\overline\Omega$. Further, let $\rho\in C^0(\overline\Omega)$ be a positive function on $\overline\Omega$ if $\gamma$ is a real valued function and complex valued function on $\overline\Omega$ with non-positive imaginary part if $\gamma$ is a complex valued function. In addition, if $\gamma$ and $\rho$ have any discontinuities in their derivatives, we assume that they are away from $\partial\Omega$. Let  $g=g(x)\in H^{3/2}(\partial\Omega)$ and $\omega^2 \geq 0$, where $H^{3/2}(\partial\Omega)$ is the Sobolev space of fractional differential order $3/2$ in the $L^2(\partial\Omega)$ sense.

\medskip
Now for given Dirichlet data $g=g(x)$, we consider the following boundary value problem:
\begin{equation}\label{eq::fp}
\left \{
\begin{array}{ll}
L_\gamma u:=\nabla\cdot\big [\gamma(x)A(x)\nabla u(x)\big ]+\omega^2\rho(x)u(x)=0 &\mbox{in } \Omega,\\
u(x)=g(x)  &\mbox{on } \partial \Omega.
\end{array}
\right.
\end{equation}

\medskip
\noindent {\it Assumption 2:} When $\gamma$ and $\rho$ are real valued function, we assume that $\Omega$ is non-vibrating, that is \eqref{eq::fp} with $g=0$ only admits a trivial solution.\\


\noindent With the above given conditions in Assumptions 1 and 2 on $\Omega,\,\gamma,\, A,\,\rho, g$ and $\omega^2$, it is well known that there exists a unique solution $u=u(x)\in H^2(\Omega)$ to \eqref{eq::fp}, where $H^2(\Omega)$ denotes the Sobolev space of differential order $2$ in the $L^2(\Omega)$ sense.\\

\noindent {\it The Inverse Problem:} We consider the following problem.  Given the {\it interior measurement} $u(x)\,(x\in\Omega)$, the  Dirichlet boundary condition $g \not\equiv 0$ on
$\partial\Omega$, the coefficients $A,\,\rho$ on $\overline\Omega$, identify $\gamma$ in $\Omega$.\\

In this paper we will mainly consider the above inverse problem whose goal is to identify $\gamma$.  However at the end of the paper we will address a somewhat easier alternate inverse problem where we assume all of the conditions above except that $\gamma$ is assumed known and our goal then will be to identify $\rho$. Natural questions, for both of our main inverse problem and the alternate inverse problem, are the uniqueness, stability and reconstruction of  $\gamma$ (or $\rho$) from the measured data.  Here our main result is a stability result for $\gamma$ given $\rho$ (or $\rho$ given $\gamma$).


There are two major backgrounds for this inverse problems. The one is coming from a newly developed imaging modality called MRE (magnetic resonance elastography) (\cite{Ehman}). It produces movies of shear waves induced by a single frequency excitation.  Once a mathematical model of these waves is determined, analysis and mathematical algorithms can be developed, based on the mathematical model. The numerical implementation of the algorithms produces diagnostically useful images with the goal of adding to noninvasive medical diagnostic capabilities.  In (\cite{Yin}), e.g., the advance, through MRE, to produce early diagnosis of liver stiffness and fibrosis is described.  MRE is in the general class of hybrid or coupled physics imaging technologies where two physical principles, here MRI and elastic wave propagation are combined to obtain a richer data set from which to obtain diagnostically rich images.

The above boundary value problem with $A$ equal to the identity matrix $I$ and $\operatorname{Im}\rho=0$ is the simplest PDE (partial differential equation) model which is used to describe a component
$u$ of a shear wave inside human tissue. In \eqref{eq::fp}, $\omega/(2\pi)$ is the angular frequency of the time harmonic vibration, $\mbox{Re}(e^{i\omega t}g(x))$, applied to a human body where $t$ is time. When $A$ is equal to the identity matrix $I$, and $\rho$ is real the functions $\rho$, $\mbox{Re}\, \gamma$ and $\mbox{Im}\, \gamma$ are the density, storage modulus and loss modulus of human tissue.

The other is coming from hydrology. This corresponds to the case $A=I$, and $\gamma$ is a real valued $C^1(\overline\Omega)$ function and $\omega=0$, then this inverse problem has been studied in hydrology.

Next we discuss some known results. In the two dimensional case Alessandrini (\cite{Alessandrini}) gave a H\"older stability result for identifying $\gamma$ by analyzing the critical set, $\{ x \in \Omega| \nabla u = 0 \}$ of, $u$, the solution to the forward problem \eqref{eq::fp}. Further, if the zero on the right hand side of the partial differential equation in \eqref{eq::fp} is replaced by a positive, H\"older continuous function, $f$, on $\overline\Omega$, Richter (\cite{Richter}) gave a Lipschitz stability result, in any dimension, for identifying $\gamma$ by showing the non-degeneracy of $\nabla u \neq 0$ and using the maximum principle. It should be noticed that the assumption, on the positivity of $f$, is a very strong assumption and is only a sufficient condition to guarantee the non-degeneracy of $\nabla u$. This assumption replaces the need to analyze the critical set. A number of other results have been established where the analysis of the critical set for a single data set is avoided by making additional hypotheses.  Marching and elliptic algorithms for recovering an unknown real coefficient, $\gamma$, when $\rho$ is real and known in the interior, $A = I$, $u$ is given in the interior, and it is propagating there to one fixed direction which means that derivative of u in this direction does not vanish there are presented in (\cite{Lin}).

Another method for addressing the problem of recovering $\gamma$, when data sets $u$ can have critical points, is the use of multiple measurements whose input can be controlled.  If we can have multiple measurements and control the input, then the reconstruction of $\gamma$ was first given by Nakamura-Jiang-Nagayasu-Cheng (\cite{NJNC}) using complex geometric optic solutions and linking them to the input data by solving a Cauchy problem.  In that paper, the regularity assumptions on $\gamma$, $\rho$ are just $\gamma\in C^2(\overline\Omega)$, $\rho\in L^\infty(\Omega)$. When it can be assumed that multiple measurements are given, a more systematic analysis, for a wide class of hybrid inverse problems, was recently done by Bal-Uhlmann (\cite{BU1}, \cite{BU2}).  In this work the given mathematical model is linearized and then: (1) a reconstruction scheme for identifying all the coefficients $\gamma, A$ and $\rho$ of the linearized operator $L_\gamma$ is presented; and (2) a Lipschitz stability result is given
when the regularity assumptions on $\gamma,\, A$ and $\rho$ are just H\"older continuous on $\overline\Omega$.

In the paper presented here, we are concerned with extending Alessandrini's result to the higher dimensional case when
$\gamma$ is a complex valued function and $\omega^2 >0$. Even in the two dimensional case and if $\gamma$ is a real valued function, the cases $\omega=0$ and $\omega^2 >0$ are quite different.  One can see that Alessandrini's argument breaks down for the case $\omega^2 >0$. As far as we know, there is not any stability result
known for the case where $\gamma$ is complex valued, where $\gamma$ can have discontinuous second derivatives and only one measurement, $u$, as opposed to multiple measurements, is given.
We will show local H\"older stability
of our inverse problem identifying $\gamma$ for the case $\gamma$ is a complex valued $C^1$ function on
$\overline{\Omega}$ and {\it piecewise analytic}
in $\Omega$ by assuming that $\rho$ is continuous on $\overline{\Omega}$ and
piecewise analytic in $\Omega$ and $A$ is positive Hermitian and analytic.

This third assumption is given more precisely as follows. We denote by $\mathcal{A}^{(k)}(\Omega)$
(resp. $\mathcal{PA}^{(k)}(\Omega)$) the set of complex valued
$C^k$ functions on $\overline{\Omega}$
which are analytic (resp. piecewise analytic) in $\Omega$ and now clarify our definition of piecewise analytic
and further assumptions on $A,\,\rho,\,\gamma$.
\\

\begin{definition} {\rm (Piecewise Analytic)} A function $f$ is piecewise analytic in $\Omega$ 
if there exists a compact subset $S$ in $\Omega$ consisting of
a finite disjoint union of closed smooth analytic hypersurfaces such that
$f$ is analytic in $ \Omega \setminus S$
and is locally extendable as an analytic function from one side of 
$S$ across to the other side. That is,
for any $x \in S$, we can find an open neighborhood $V$ of $x$ with $V \setminus S$
consisting of two connected components $V_1$ and $V_2$ for which $f$
in $V_1$ (resp. $V_2$) analytically extends to $V$.
\end{definition}

\medskip
\noindent
{\it Assumption 3:} Let $A\in M(n;\,\mathcal{A}^{(1)}(\Omega))$, $\gamma\in\mathcal{PA}^{(1)}(\Omega)$, $\rho\in\mathcal{PA}^{(0)}(\Omega)$ and the locations of the singularities of $\gamma$ and $\rho$ be the same. Here $A\in M(n;\,\mathcal{A}^{(1)}(\Omega))$ means that all the entries of $n\times n$ matrix $A$ belong to $\mathcal{A}^{(1)}(\Omega)$.

\medskip
With {\it Assumptions 1,2, and 3} it can be shown by using the theories of analytic pseudo-differential operators, the theories for coercive boundary value problems and the fact that when $\gamma$ and $\rho$ are real we assume that we have
a nonvibrating problem, that the unique solution $u\in H^1(\Omega)$ to
\eqref{eq::fp} belongs to $\mathcal{PA}^{(1)}(\Omega)$.
This follows because the interior transmission problem can be transformed to a coercive boundary value problem for a system of equations by introducing the boundary normal coordinates in the neighborhood $V$ of $x\in S$ (see Definition above) so that 
we can reflect the component $V_1 \subset V$ to the other side of $S$ where we have the component $V_2 \subset V$, and then apply the analytic hypo-ellipticity result given in Chapters III and V of \cite{Treves} to this coercive boundary value problem.

\medskip
Now let $\sigma > 0$ be a sufficiently small constant.  Furthermore, we introduce the notion of an admissible pair for functions on $\overline{\Omega}$.\\

\begin{definition} {\rm (Admissible Pair)}
A pair $(\gamma_1,\, \gamma_2)$ of functions on $\overline{\Omega}$
is said to be admissible if  there exist exceptional angles
$
\kappa_1 < \kappa_2 < \dots < \kappa_\ell < \kappa_{\ell + 1} = \kappa_1 + 2\pi
$
with $\kappa_{k+1} - \kappa_k \le \pi - \sigma$ $(k=1, \dots, \ell)$ such that,
for any $k=1,\dots,\ell$,
\begin{equation}{\label{eq:def-addmissible}}
H_{n-1}\left(\left\{x \in \overline{\Omega}\,:\,
(\gamma_2 - \gamma_1)(x) \ne 0,\,\,
\operatorname{arg}\, (\gamma_2 - \gamma_1)(x) \equiv \kappa_k
\mod 2\pi\right\}\right) = 0.
\end{equation}
Here $H_{n-1}$ denotes the $(n-1)$-dimensional Hausdorff measure.
\end{definition}


\

Clearly the fact that  $\ell \ge 3$ follows from definition.  In addition, we make the following remarks.

\medskip
\noindent
{\it{Remark}:}
\newline
(i) Here we give a sufficient condition in order that $\gamma_j, \, j\,=\,1,2$ is an admissible pair:
Suppose there exists a non-negative constant $\kappa < (\pi - \sigma)/2$ satisfying
\begin{equation}
|\operatorname{Im}\,(\gamma_2 - \gamma_1)(x)|
\le \tan(\kappa)\, |\operatorname{Re}\,(\gamma_2 - \gamma_1)(x)|
\qquad (x \in \overline{\Omega}),
\end{equation}
then the pair $(\gamma_1,\, \gamma_2)$
becomes admissible.  Note also, as a particular case,  real valued $(\gamma_1,\, \gamma_2)$ are always
admissible. Furthermore $\gamma_1, \gamma_2$ are always admissible if $\operatorname{Im}(\gamma_2 - \gamma_1) \equiv 0$.
\newline
(ii) In addition to the lower bound on $\ell$, given in the definition of Admissible Pair, we can also assume without loosing its generality that $\ell \leq 4$.  To see this suppose that $\kappa_1,\dots,\kappa_\ell$ are the exceptional angles of the
admissible pair $(\gamma_1,\,\gamma_2)$
for which (\ref{eq:def-addmissible}) holds.
We will show that if $\ell \ge 5$, there exists a
$k$ with $\kappa_{k+2} - \kappa_k \le \pi - \sigma$ so that $\kappa_{k+1}$ can be eliminated.
Suppose that no such $k$ exists for some $\ell \ge 5$.
Then $\kappa_{k+2} - \kappa_k > \pi - \sigma$ holds for every $k$ where we set
$\kappa_{k+\ell} := \kappa_k + 2\pi$ for convenience. Clearly we have
$$
2(\kappa_{\ell + 1} - \kappa_1) = \sum_{k=1}^\ell(\kappa_{k+2} - \kappa_k)
> \ell (\pi - \sigma) \ge 5 (\pi - \sigma),
$$
which contradicts the fact $\kappa_{\ell+1} - \kappa_1 = 2\pi$.

\

In order to state our main result, it is convenient to use the Sobolev space ${W^{q,p}(\Omega)}$ of differential order $q$ in the $L^p(\Omega)$ sense with the norm $||\cdot||_{W^{q,p}(\Omega)}$. We also set $||B||_{W^{q,p}(\Omega)} := \sum_{i,j}
||b_{i,j}||_{W^{q,p}(\Omega)}$ for a matrix $B(x) = (b_{i,j}(x))$.  We denote by $u_k\in H^2(\Omega)\,(k=1,2)$ the solutions to the boundary value problem
\eqref{eq::fp} with $\gamma=\gamma_k,\,g=g_k\not\equiv0\,(k=1,2)$ and $L_k=L_{\gamma_k}\,(k=1,2)$.


\medskip
Then, we have our main result.

\medskip
\begin{theorem}{\rm (Main Theorem)}
{\it Let $d > 0$ and $\Omega_d := \{x \in \Omega\,:\, \operatorname{dist}(x,\, \mathbb{R}^n \setminus \Omega) > d\}$. Let $\gamma_1, \gamma_2$ be an admissible pair and $A, \rho, \gamma_1, \gamma_2, \omega^2$  satisfy Assumptions 1,2,3, with $||A||_{W^{1,\infty}(\Omega)} \le \sigma^{-1}$ and
$||\gamma_k||_{W^{1,\infty}(\Omega)} \le \sigma^{-1},\, k=1,2$.
Then there exist constants $C>0$ and $\alpha\in(0,1)$ depending only on $\Omega$, $d$,
$\sigma$, $g_1$ and the coefficients of $L_1$ such that
\begin{equation}\label{eq::stability}
\Vert\gamma_2-\gamma_1\Vert_{L^\infty(\Omega_d)}
\le C\left(\Vert\gamma_2-\gamma_1\Vert_{L^\infty(\partial\Omega)}
+\Vert u_2-u_1\Vert_{W^{2,1}(\Omega)}\right)^\alpha
\end{equation}
for any $g_2$ and $\gamma_2$ of an admissible pair $(\gamma_1,\, \gamma_2)$.
Furthermore, if $\Omega$ has an analytic smooth boundary and if
$g_1$ is analytic in $\partial \Omega$ and all the coefficients of $L_1$ are
analytic near $\partial \Omega$,
then we have the estimate} (\ref{eq::stability}) {\emph in which $L^\infty(\Omega_d)$ is replaced with
$L^\infty(\Omega)$.}\\
\end{theorem}

The succeeding sections are devoted to the proof of the main result and they are organized as follows.
We first present a key identity and an associated estimate. Then,
we give: (1) statements about a tubular neighborhood of the critical set of the solution, $u_1$, to \eqref{eq::fp}, when $\gamma$ is replaced by $\gamma_1$; (2) the estimate of the $n$-dimensional Lebesgue measure of the tubular neighborhood; and (3) a lower estimate of $|\nabla u|$ outside this tubular neighborhood. The proofs of these results are given in Appendix. Combining the three sets of estimates we finish proving the main result. Finally we give the stability estimate for the alternate inverse problem which is to identify $\rho$ given $\gamma, A, \omega^2$, and $u$.

\section{Key identity and an associated estimate}
Let $\psi=\gamma_2-\gamma_1$ and $\nabla_{B}$ denote $B(x) \nabla$ for a matrix $B(x)$.
Then, it is straightforward to establish the
following key identity.\\

\begin{lemma}{\rm (Key Identity)}
Let Assumptions 1 and 2 be satisfied with $\gamma$ replaced by $\gamma_k$ and $g$ replaced by $g_k, \, k=1,2$. Then, for any $\zeta\in H_0^1(\Omega)$, that is $\zeta\in H^1(\Omega)$ with trace $\zeta\Big|_{\partial\Omega}=0$, we have
\begin{equation}\label{eq::identity}
	\int_\Omega\psi\Anabla u_1\cdot\CAnabla \zeta=-\int_\Omega
	\gamma_2\Anabla (u_2-u_1)\cdot\CAnabla \zeta+
\omega^2\int_\Omega\rho(u_2-u_1)\zeta.
\end{equation}
Here $x \cdot y$
denotes a sum $\sum_{k=1}^n x_ky_k$ for $x,y \in \mathbb{C}^n$.
\end{lemma}

%

\vspace{.1in}
Based on this key identity, we have the following fundamental estimate associated to the key identity.\\

\begin{proposition}\label{asso_estimate}
Let Assumptions 1,2, and 3 be satisfied with $\gamma$ replaced by $\gamma_k$ and $g$ replaced by $g_k,\,k=1,2$. Then, there exists a constant $C>0$ depending only on $\Omega$,  $\sigma$,
$g_1$ and the coefficients of $L_1$ such that
\begin{equation}\label{important_estimate}
	\int_\Omega |\psi|\,\Unabla \le
C\left(\Vert\psi\Vert_{L^\infty(\partial\Omega)}+\Vert u_2-u_1\Vert_{W^{2,1}(\Omega)}
\right).
\end{equation}
Note that, as $A$ is a positive Hermitian matrix,
the term $\Unabla$ takes non-negative real values.
\end{proposition}

\medskip
Before we present the proof, taking advantage of our Assumptions 1 and 3, we first present several new sets and functions that we need for the estimate.
Consider the map $\iota: \partial \Omega \times \mathbb{R} \to \mathbb{R}^n$
of $C^1$ class defined by
$(y,s) \to y + s \nu(y)$ where $\nu(y)$ is a unit conormal vector of
$\partial \Omega$ at $y$ pointing to $\Omega$. Then, as $\partial \Omega$ is compact,
there exists an $\epsilon > 0$ such that $\iota$ becomes a $C^1$ isomorphism between
$\partial \Omega \times (-\epsilon, \epsilon)$ and
an open neighborhood of $\partial \Omega$. Hence we have a family
$\{U_j\}_{j \in \mathbb{N}}$ of relatively compact open subsets in $\Omega$
satisfying the conditions below:
\begin{enumerate}
\item $\Omega = \cup_j U_j$ and $U_j \subset\subset U_{j+1} \subset\subset \Omega$.
\item $U_j$ has a $C^1$ smooth boundary.
\item $H_{n-1}(\partial U_j) \to H_{n-1}(\partial \Omega)$ ($j \to \infty$). We also have
$\operatorname{dist}(\partial \Omega,\, \partial U_j) \to 0$ when $j \to \infty$.
\end{enumerate}
Furthermore, using the definition of subanalytic sets given in Appendix (see also (\cite{ep}), we can find a family $\{W_j\}_{j \in \mathbb{N}}$ of
relatively compact subanalytic open subsets in $\Omega$ with
$U_j \subset\subset W_j \subset\subset \Omega$.
One choice for $W_j$ can be obtained by dividing $\mathbb{R}^n$ into sufficiently small n-dimensional cubes. Then $W_j$ can be selected to be a finite union of these cubes where each cube intersects $\overline{U_j}$ and
where the closure of the finite union is contained in $\Omega$.


Now divide the complex plane, $\mathbb{C}$, into proper sectors
\begin{equation}
\Gamma_k := \{z\in\mathbb{C}\setminus\{0\}\,:\,
\kappa_k \le \mbox{arg}\, z \le \kappa_{k+1}\} \cup \{0\}
\qquad (k=1,\dots,\ell).
\end{equation}
Set
\begin{equation}
\widehat{\Omega}_k := \{x \in \overline{\Omega}\,:\, \psi(x) \in \Gamma_k \},
\end{equation}
and let $\theta_k(x)$ be the Lipschitz continuous function
\begin{equation}
\theta_k(x) :=
\left\{
\begin{array}{ll}
\operatorname{Re}\,\psi_k(x) - c_k\operatorname{Im}\,\psi_k(x) \qquad
& \text{if } \operatorname{Im}\,\psi_k(x) \ge 0, \\
\\
\operatorname{Re}\,\psi_k(x) + c_k\operatorname{Im}\,\psi_k(x) \qquad
& \text{if } \operatorname{Im}\,\psi_k(x) \le 0,
\end{array}
\right.
\end{equation}
where $\psi_k(x) := \exp(-(\kappa_k + \kappa_{k+1})i/2)\psi(x)$ and
$c_k := 1/\tan( (\kappa_{k+1} - \kappa_k)/2)$.
In addition, using the definition of subanalytic functions in Appendix, see also (\cite{ep}), $\theta_k$
is a subanalytic function on $\overline{W_j}$.

Furthermore it follows from the definition of $\theta_k(x)$ that we have
\begin{enumerate}
\item $\theta_k(x) < 0$ if and only if $\psi(x) \in \mathbb{C} \setminus \Gamma_k$,
\item $\theta_k(x) = 0$ if and only if $\psi(x) \in \partial \Gamma_k$,
\item $\theta_k(x) > 0$ if and only if $\psi(x) \in \Gamma_k^\circ$.
\end{enumerate}
Hence, in particular, a point $x$ belongs to $\widehat{\Omega}_k$ if and only
if $\theta_k(x) \ge 0$ holds.
We also define, for $h>0$,
\begin{equation}
\theta_{k,h}(x) := h^{-1}([\theta_k(x)]^+\wedge h),
\end{equation}
where $[m]^+=\mbox{max}(m,\,0)$, $m\wedge\ell=\mbox{min}(m,\,\ell)$
for $m,\,\ell\in\mathbb{R}$. Then $\theta_{k,h}(x)$ is again
a Lipschitz continuous subanalytic function on each $\overline{W}_j$ and
it satisfies
\begin{equation}\label{eq:base_prop_of_theta_kh}
0 \le \theta_{k,h}(x) \le 1,\,\,\,
\operatorname{supp} \theta_{k,h}(x) \subset \widehat{\Omega}_k,\,\,\,
\lim_{h \to 0^+} \theta_{k,h}(x) \to
\chi_{\{x \in \widehat{\Omega}_k\,:\, \theta_k(x) \ne 0\}}.
\end{equation}
Here $\chi_A$ designates the characteristic
function of a subset $A$.

Finally, let $\tau_{j,h}$ be a $C^\infty$ function in $\Omega$ satisfying
$0 \le \tau_{j,h}(x) \le 1$,
$\tau_{j,h}(x) = 0$ at $x \in \Omega$ with
$\operatorname{dist}(x,\, \mathbb{R}^n \setminus U_j) \le 2^{-1}h$,
$\tau_{j,h}(x) = 1$ at $x \in \Omega$ with
$\operatorname{dist}(x,\, \mathbb{R}^n \setminus U_j) \ge h$
and $|\nabla \tau_{j,h}(x)| \le C_\tau h^{-1}$ for some $C_\tau > 0$.
Note that, by choosing a suitable $\tau_{j,h}$ for each $j$,
we may assume that the constant $C_\tau$ is independent of $j$.

\medskip
Now we are ready to prove Proposition \ref{asso_estimate}.  In the proof we will make extensive use of the
results in Appendix.\\

\begin{proof}
Set
\begin{equation}
	\zeta_{j,k,h} := \overline{u}_1 \tau_{j,h} \theta_{k,h} \qquad (k=1,\dots,\ell,\,j \in \mathbb{N}).
\end{equation}
As $\zeta_{j,k,h}$
belongs to $H_0^1(\Omega)$ and $\operatorname{supp} \zeta_{j,k,h} \subset U_j$,
by taking $\zeta_{j,k,h}$ as $\zeta$ in the key identity, and replacing $\Omega$ by $U_j$,
we obtain
\begin{equation} \label{eq:fundamental equality}
\begin{aligned}
	\left|\,\int_{U_j} \psi \tau_{j,h} \theta_{k,h} \Unabla \right|
&\le
\left|\,\int_{U_j} \psi \overline{u}_1 \tau_{j,h} \Anabla u_1 \cdot \CAnabla \theta_{k,h} \right|
+ \left|\,\int_{U_j} \psi \overline{u}_1 \theta_{k,h}
\Anabla u_1 \cdot \CAnabla \tau_{j,h} \right| \\
&\quad
+ \left|\,\int_{U_j} \gamma_2
\Anabla (u_2 - u_1) \cdot \CAnabla \zeta_{j,k,h} \right|
+ \omega^2 \left|\,\int_{U_j} \rho (u_2 - u_1)\zeta_{j,k,h} \right|.
\end{aligned}
\end{equation}

We will compute, when $h \to 0^+$,
the limit of each of the five terms in (\ref{eq:fundamental equality}), or the limit of estimates of each of the five terms.\\

\noindent {\it Case 1: The limit of the term on the left hand side of (\ref{eq:fundamental equality}).}\\

 \noindent It follows from (\ref{eq:base_prop_of_theta_kh}) that
we have
$$
\lim_{h\rightarrow 0^+}\int_{U_j} \psi \tau_{j,h} \theta_{k,h} \Unabla
= \int_{\{x \in \widehat{\Omega}_k \cap U_j\,:\, \theta_k(x) \ne 0\}}
\psi \Unabla
$$
Set $T := \left\{x \in \overline{W}_j\,:\, \theta_k(x) = 0,\,
\psi(x) \ne 0\right\}$ which is a relatively compact subanalytic subset in $\mathbb{R}^n$.
Then,
by the admissible condition for the pair $(\gamma_1,\, \gamma_2)$,
we have
$
\dim_{\mathbb{R}} T < n-1.
$
Hence the $n$-dimensional volume of
$T$ is zero (see Proposition \ref{prop:area_finite} in Appendix)
and we can conclude
$$
\lim_{h\rightarrow 0^+}\int_{U_j} \psi \tau_{j,h} \theta_{k,h} \Unabla
= \int_{\{x \in \widehat{\Omega}_k \cap U_j\,:\, \theta_k(x) \ne 0\}}
\psi \Unabla
=
\int_{\widehat{\Omega}_k \cap U_j} \psi \Unabla
$$

\noindent {\it Case 2: The limit of the first term in the right-hand side of (\ref{eq:fundamental equality}).}\\

\noindent  We will establish that this term tends
to zero when $h \to 0^+$. To obtain our result we first note that
it follows from Proposition \ref{prop:area_finite}
in Appendix that there
exists a constant $M_{\theta_k}$ such that, for any $t \in \mathbb{R}$ with
finite exceptional ones, we get
$$
H_{n-1}(\{x \in U_j\,:\, \theta_k(x) = t\}) \le
H_{n-1}(\{x \in W_j\,:\, \theta_k(x) = t\}) \le M_{\theta_k}.
$$
Set
$$
U_{j,k,h} := \{x \in U_j\,:\, 0 < \theta_k(x) < h\}.
$$
Then, by noticing the fact that $\theta_{k,h}(x)$ is locally constant
in $U_j \setminus \overline{U_{j,k,h}}$, we have
$$
\begin{aligned}
\left|\,\,\int_{U_j} \psi \overline{u}_1 \tau_{j,h} \Anabla u_1 \cdot \CAnabla
\theta_{k,h} \right|
&=
\left|\,\,\int_{\overline{U_{j,k,h}} \cap U_j}
\psi \overline{u}_1 \tau_{j,h} \Anabla u_1 \cdot \CAnabla \theta_{k,h} \right| \\
&=
\left|\,\int_{U_{j,k,h}}
\psi \overline{u}_1 \tau_{j,h} \Anabla u_1 \cdot \CAnabla \theta_{k,h} \right|.
\end{aligned}
$$
Here the last equality follows from the fact that
$\partial U_{j,k,h} \cap U_j$ has n-dimensional measure zero.
To see this let $W_{j,k,h}:= \{x \in W_j\,:\, 0 < \theta_k(x) < h\}$. Then
as $W_{j,k,h}$ is a subanalytic subset, we have
$\operatorname{dim}_{\mathbb{R}} \partial W_{j,k,h} < n$.
Hence the $n$-dimensional volume of $\partial W_{j,k,h}$
is zero.  This can be seen, for example, by again utilizing Proposition \ref{prop:area_finite} in Appendix. Then the inclusion
$\partial U_{j,k,h} \cap U_j \subset \partial W_{j,k,h} \cap U_j$ yields the result.

Set 
$$
Z^{0} := \{x \in W_j\,:\, \theta_k(x) = 0,\, \psi(x) = 0\},\,\,
Z^{\times} := \overline{\{x \in W_j\,:\, \theta_k(x) = 0,\, \psi(x) \ne 0\}}
$$
and we also set, for $\epsilon > 0$,
$$
Z_\epsilon^{0} := \{x \in W_j\,:\, \operatorname{dist}(x,Z^{0}) \le \epsilon\},\,\, 
Z_\epsilon^{\times} := \{x \in W_j\,:\,
\operatorname{dist}(x,Z^{\times}) \le \epsilon\}.
$$
Note that, by the admissible pair condition, we have
$\operatorname{dim}_{\mathbb{R}}(Z^{\times})
< n -1$.
In the estimate
$$
\begin{aligned}
\left|\,\int_{U_{j,k,h}}
\psi \overline{u}_1 \tau_{j,h} \Anabla u_1 \cdot \CAnabla \theta_{k,h} \right|
&\le
\left|\,\int_{U_{j,k,h} \cap Z^{0}_\epsilon}
\psi \overline{u}_1 \tau_{j,h} \Anabla u_1 \cdot \CAnabla \theta_{k,h} \right| \\
&\qquad+
\left|\,\int_{U_{j,k,h} \setminus Z^{0}_\epsilon}
\psi \overline{u}_1 \tau_{j,h} \Anabla u_1 \cdot \CAnabla \theta_{k,h} \right|,
\end{aligned}
$$
by using the co-area formula,
we obtain estimates
\begin{equation}
\begin{aligned}
	\left|\,\int_{U_{j,h,k} \cap Z^{0}_\epsilon} \psi \overline{u}_1 \tau_{j,h} \Anabla u_1 \cdot \CAnabla
\theta_{k,h} \right|
&\le
\left(\sup_{x \in Z^0_\epsilon}|\psi(x)|\right)
||(\overline{u}_1|\Anabla u_1|)||_{L^\infty(\Omega)}
\dfrac{1}{h}\int_{U_{j,k,h}} |\CAnabla \theta_k| \\
&\le
\left(\sup_{x \in Z^0_\epsilon}|\psi(x)|\right)
||(\overline{u}_1|\Anabla u_1|)||_{L^\infty(\Omega)} M_{\theta_k}
\end{aligned}
\end{equation}
and
\begin{equation}\label{eq:estimate-boundary-non-zero}
\begin{aligned}
&\left|\,\int_{U_{j,k,h} \setminus Z^0_\epsilon} \psi \overline{u}_1 \tau_{j,h} \Anabla u_1 \cdot \CAnabla
\theta_{k,h} \right|
\le
2\sigma^{-1}
||(\overline{u}_1|\Anabla u_1|)||_{L^\infty(\Omega)}
\dfrac{1}{h}\int_{U_{j,k,h} \setminus Z^0_\epsilon} |\CAnabla \theta_k| \\
&\quad\le
2\sigma^{-1}
||(\overline{u}_1|\Anabla u_1|)||_{L^\infty(\Omega)}
\left|\left|H_{n-1}\left(\left\{x \in U_{j} \setminus Z^0_\epsilon\,:\, \theta_k(x) = t
\right\}\right)\right|\right|_{L^\infty(\{t \in \mathbb{R}\,:\,0 < t < h\})}.
\end{aligned}
\end{equation}
For any $\epsilon' > 0$, we have
$U_{j,k,h} \setminus Z^0_\epsilon \subset Z^\times_{\epsilon'} \setminus Z^0_\epsilon$ if
$h > 0$ is sufficiently small. Therefore, by the fact that
$\operatorname{dim}_{\mathbb{R}}Z^\times < n-1$,
it follows from the second claim of
Proposition \ref{prop:area_finite}
in Appendix that the right-hand side of
(\ref{eq:estimate-boundary-non-zero}) tends to zero when $h \to 0^+$.
Hence we have obtained
$$
\overline{\underset{h \to 0^+}{\lim}}
\left|\,\int_{U_j} \psi \overline{u}_1 \tau_{j,h} \Anabla u_1 \cdot \CAnabla
\theta_{k,h} \right|
\le
\left(\sup_{x \in Z^0_\epsilon}|\psi(x)|\right)
||(\overline{u}_1|\Anabla u_1|)||_{L^\infty(\Omega)} M_{\theta_k}.
$$
Clearly we have
$
\displaystyle\sup_{x \in Z^0_\epsilon}|\psi(x)| \to 0
$
$(\epsilon \to 0^+)$, from which
$$
\lim_{h \to 0^+} \left|\,\int_{U_j} \psi \overline{u}_1 \tau_{j,h} \Anabla u_1 \cdot \CAnabla \theta_{k,h} \right| = 0 \,\,\,\;\;\;\mbox{immediately follows}.
$$
{\it Case 3: The limit of the second term on the right hand side of (\ref{eq:fundamental equality}).}\\

\noindent Set
$$
\overline{U}_{j,h} :=
\{x \in \overline{U_j}\,:\,
\operatorname{dist}(x,\,\mathbb{R}^n \setminus U_j) \le h\}.
$$
We have
$$
\begin{aligned}
	\left|\,\int_{U_j} \psi \overline{u}_1 \theta_{k,h}
	\Anabla u_1 \cdot \CAnabla \tau_{j,h} \right|
& =
\left|\,\,\int_{\overline{U}_{j,h} \cap U_j} \psi \overline{u}_1 \theta_{k,h}
\Anabla u_1 \cdot \CAnabla \tau_{j,h} \right|  \\
& \le
\left(\sup_{x \in \overline{U}_{j,h}} |\psi(x)| \right)
||(\overline{u}_1 |\Anabla u_1|)||_{L^\infty(\Omega)}
\dfrac{C_\tau}{h} \operatorname{vol}(\overline{U}_{j,h}).
\end{aligned}
$$
As $\partial U_j$ is $C^1$ smooth,
we get
$\underset{h \to 0^+}{\lim}
h^{-1}\operatorname{vol}(\overline{U}_{j,h}) = H_{n-1}(\partial U_j)$.
Hence we have
$$
\underset{h \to 0^+}{\overline{\lim}}
\left|\,\int_{U_j} \psi \overline{u}_1 \theta_{k,h}
\Anabla u_1 \cdot \CAnabla \tau_{j,h} \right|
\le
\left(C_\tau H_{n-1}(\partial U_j)
||(\overline{u}_1 |\Anabla u_1|)||_{L^\infty(\Omega)}\right)
||\psi||_{L^\infty(\partial U_j)}.
$$

\noindent {\it {Case 4:  The limit of the last two terms of}} (\ref{eq:fundamental equality}).\\

\noindent Our Assumptions imply that $u_2 - u_1$ is in $H^2(\Omega)=W^{2,2}(\Omega)$. Hence, we have
$$
\begin{aligned}
&\quad\left|\,\,\int_{U_j} \gamma_2
\Anabla (u_2 - u_1) \cdot \CAnabla \zeta_{j,k,h} \right|
+ \omega^2 \left|\,\,\int_{U_j} \rho (u_2 - u_1)\zeta_{j,k,h} \right| \\
&\qquad\qquad=
\left|\,\,\int_{U_j}
\operatorname{div}(\gamma_2
\nabla_A (u_2 - u_1)) \zeta_{j,k,h} \right|
+ \omega^2 \left|\,\,\int_{U_j} \rho (u_2 - u_1)\zeta_{j,k,h} \right| \\
&\qquad\qquad\le
\left(2\sigma^{-2} + \omega^2||\rho||_{L^\infty(\Omega)}\right)
||u_2 - u_1||_{{W^{2,1}}(\Omega)}.
\end{aligned}
$$

\noindent {\it The last step: Combining Cases 1-4}. \\

\noindent Summing up, we have, for $k=1, \dots, \ell$,
\begin{equation}
\left|\,\, \int_{\widehat{\Omega}_k \cap U_j} \psi \Unabla\,\, \right|
\le
C_j \left(|| \psi ||_{L^\infty(\partial U_j)} +
||u_2 - u_1||_{W^{2,1}(\Omega)}\right),
\end{equation}
where
$
C_j := \max
\left\{C_\tau H_{n-1}(\partial U_j)
||(\overline{u}_1 |\Anabla u_1|)||_{L^\infty(\Omega)},\,\,
2\sigma^{-2} + \omega^2||\rho||_{L^\infty(\Omega)}\right\}.
$
Now, by letting $j \to \infty$, we obtain
\begin{equation}
\left|\,\, \int_{\widehat{\Omega}_k} \psi \Unabla\,\, \right|
\le
C \left(|| \psi ||_{L^\infty(\partial \Omega)} +
||u_2 - u_1||_{W^{2,1}(\Omega)}\right)
\end{equation}
with
\begin{equation}
C := \max
\left\{C_\tau H_{n-1}(\partial \Omega)
||(\overline{u}_1 |\Anabla u_1|)||_{L^\infty(\Omega)},\,\,
2\sigma^{-2} + \omega^2||\rho||_{L^\infty(\Omega)}\right\}.
\end{equation}
Then, by noticing the fact
\begin{equation}
	|\psi(x)| \le \dfrac{1}{\sin(\sigma/2)}
	\operatorname{Re}\, \beta_k \psi(x)\quad
(x \in \widehat{\Omega}_k)
\end{equation}
with $\beta_k := \exp(-(\kappa_{k+1}+\kappa_k)i/2)$
because of $\psi(x) \in \Gamma_k$ for any $x \in \widehat{\Omega}_k$,
we obtain
$$
\begin{aligned}
&\int_\Omega |\psi| \Unabla
\le  \sum_{k=1}^\ell
\int_{\widehat{\Omega}_k} |\psi| \Unabla
\le \dfrac{1}{\sin(\sigma/2)} \sum_{k=1}^\ell \operatorname{Re}
\int_{\widehat{\Omega}_k} \beta_k \psi \Unabla \\
&\quad= \dfrac{1}{\sin(\sigma/2)} \sum_{k=1}^\ell \left | \operatorname{Re}
 \int_{\widehat{\Omega}_k} \beta_k \psi \Unabla \right|
 \le \dfrac{1}{\sin(\sigma/2)} \sum_{k=1}^\ell \left |\,
 \int_{\widehat{\Omega}_k} \beta_k \psi \Unabla \right| \\
&\quad= \dfrac{1}{\sin(\sigma/2)} \sum_{k=1}^\ell \left |
\,\int_{\widehat{\Omega}_k} \psi \Unabla \right|
 \le \dfrac{4C}{\sin(\sigma/2)} \left(
|| \psi ||_{L^\infty(\partial \Omega)} +
||u_2 - u_1||_{W^{2,1}(\Omega)}\right).
\end{aligned}
$$
This completes the proof. 
\end{proof}

\section{The Critical set of  $u_1$}
We consider first the case where $\Omega$ has $C^2$ boundary $\partial \Omega$.  In this case
let $V$ and $W$ be relatively compact open subsets in $\Omega$ satisfying
\begin{enumerate}
\item $\Omega_d \subset V \subset W \subset\subset \Omega$.
\item $V$ has a $C^1$ smooth boundary.
\item $W$ is a subanalytic subset in $\mathbb{R}^n$.
\end{enumerate}
These $V$ and $W$ can be constructed by using the argument following the statements before the proof of Proposition \ref{asso_estimate}.\\

\noindent
\begin{lemma}{\rm (The Tubular Neighborhood of the Critical Set)}
Let $Z$ be the critical set of $u_1$ in $\overline{W}$. That is
\begin{equation}
Z=\{x\in\overline{W}\,:\, \nabla u_1(x)=0\}.
\end{equation}
Then, there exist a family $U(\eta)\,(0<\eta\le 1)$ of subanalytic open neighborhoods
of $Z$, positive constants $r$ and $C_j\,(1\le j\le 3)$, where $C_j$ is
independent of $\eta$ such that
\begin{itemize}
\item[1.]$\operatorname{vol}(U(\eta))\le C_1\eta\quad\,\,(\eta\in(0,1]).$
\item[2.] $\operatorname{dist}(x,\,Z)\ge C_2\eta\quad\,\,(x\in \mathbb{R}^n\setminus
U(\eta),\,\eta\in(0,1]).$
\item[3.] $\Unabla \ge C_3\operatorname{dist}(x,\,Z)^r\quad\,\,(x\in\overline{V}).$
\end{itemize}
\end{lemma}

\medskip
\begin{proof}
We first note that, as $u_1$ is piecewise analytic in an open neighborhood
of $\overline{W}$, the function
$\Unabla$ is subanalytic on $\overline{W}$ and
$Z$ is a compact subanalytic subset in $\mathbb{R}^n$.
Furthermore, since $g_1$ is non-zero,
by the unique continuation property of a solution for $L_1$,
we have $\operatorname{dim}_{\mathbb{R}} Z < n$.
Hence, by the two Theorems
in Appendix,
for every $\eta,\; 0<\eta\le 1$,
there exists $U(\eta)$, a subanalytic open neighborhood
of $Z$ and positive constants $r$ and $C_j\,(1\le j\le 3)$ such that
\begin{itemize}
\item[1.]$\mbox{vol}(U(\eta))\le C_1\eta\quad\,\,(\eta\in(0,1]).$
\item[2.] $\mbox{dist}(x,\,Z)\ge C_2\eta\quad\,\,(x\in \mathbb{R}^n\setminus
U(\eta),\,\eta\in(0,1]).$
\item[3.] $\Unabla \ge C_3\mbox{dist}(x,\,Z)^r\quad\,\,(x\in\overline{V})$.
\end{itemize}
Hence the proof is complete.     \end{proof}

\medskip
When $\partial \Omega$ is analytic smooth, $g_1$ is analytic in $\partial \Omega$ and all the coefficiens of $L_1$ are analytic near $\partial \Omega$, $\Unabla$ is a subanalytic function on $\overline{\Omega}$ and the subset $Z=\{x\in\overline{W}\,:\, \nabla u_1(x)=0\}$ is compact and subanalytic in ${\mathbb{R}}^n$.  In this case then the conclusions of the
lemma hold in all of $\overline{\Omega}$.\\

\section{The final steps in the proof of the Main Theorem}

We will use the estimate in Proposition \ref{asso_estimate} and the properties of the critical set of $u_1$, given above,
to estimate $\int_{\Omega_d}|\psi|$ when $\partial \Omega$ is a $C^2$ boundary of $\Omega$ and to estimate $\int_{\Omega}|\psi|$ when
$\partial \Omega$ is of $C^2$.  We begin with the case where $\partial \Omega \in C^2$ and let $V$ be as
defined in the previous section.
In this case, for any $\eta\in(0,1]$, we have
\begin{equation}\label{eq::L1estimate1}
\begin{array}{ll}
\int_{V}|\psi|&\le\int_{V \cap U(\eta)}|\psi|+\int_{V \setminus U(\eta)}|\psi|\\
\\
&\le  C_1\eta\Vert\psi\Vert_{L^\infty(V)}
+( C_2^{r}C_3)^{-1}\eta^{-r}\int_{V}|\psi|\,\Unabla.
\end{array}
\end{equation}
Then, minimizing
\eqref{eq::L1estimate1} with respect to $\eta\in(0,1]$ and possibly making $C_1$ larger in order to ensure that $\eta\;\in\;(0,1]$
we have
\begin{equation}\label{eq::L1estimate2}
\int_{V}|\psi|\le C\Vert\psi\Vert_{L^\infty(V)}^{r/(r+1)}
\left(\int_V|\psi|\,\Unabla\right)^{1/(r+1)}
\end{equation}
for some constant $C>0$ depending only on $g_1$, $V$, and the coefficients of $L_1$.

As $V$ satisfies the cone condition, by the Gagliardo-Nirenberg inequality, there exists
a constant $C\rq{}>0$ depending only on $V$ such that for any $s>n$, we have
\begin{equation}\label{eq::GNineq}
\Vert\psi\Vert_{L^\infty(V)}\le
C\rq{}\Vert\psi\Vert_{W^{1,s}(V)}^\theta
\Vert\psi\Vert_{L^1(V)}^{1-\theta}\le {C_\psi}\rq{}\Vert\psi
\Vert_{L^1(V)}^\kappa,
\end{equation}
where $\theta=n/(n+1-n/s)$, ${C_\psi}\rq{}=C\rq{}\Vert\psi
\Vert_{W^{1,s}(V)}^\theta$ and $\kappa=1-\theta$.
Combining \eqref{eq::L1estimate2} and \eqref{eq::GNineq},
we have
\begin{equation}
\begin{aligned}
\Vert\psi\Vert_{L^\infty(\Omega_d)}\le
\Vert\psi\Vert_{L^\infty(V)}
&\le
C_\psi
\left(\int_V|\psi|\, \Unabla \right)^{\kappa/(r+1)} \\
&\qquad\le
C_\psi
\left(\int_\Omega |\psi|\, \Unabla \right)^{\kappa/(r+1)},
\end{aligned}
\end{equation}
where $C_\psi=C_\psi\rq{}\left(C ||\psi||^{r/(r+1)}_{L^\infty(V)}\right)^{\kappa}$.
Therefore the estimate \eqref{eq::stability} immediately follows from 
Proposition \ref{asso_estimate}.

Finally we show the last assertion of the theorem. Since the solution $u_1$ becomes,
in this case, piecewise analytic in an open neighborhood of $\overline\Omega$ and since
$\Omega$ itself is subanalytic,
$\Unabla$ is a subanalytic function on $\overline{\Omega}$ and
the subset $Z:=\{x\in\overline{\Omega}\,:\, \nabla u_1(x)=0\}$
is compact and subanalytic in $\mathbb{R}^n$. Hence the same argument in this section
can be applied to the case $V=W=\Omega$, and we have obtained the final estimate with $V=\Omega_d=\Omega$.
In this case the exponent on the right hand side can be left the same or also changed to $(1-\theta)/(1+r\theta)$ with
$C_\psi=C^\prime_\psi(C)^\kappa$.

\section{Local stability for $\rho$ given $\gamma$}

In this section we will consider the alternate inverse problem as stated in the introduction. That is we consider the inverse problem of identifying $\rho$, given $\gamma,\;A,\;\omega^2$, an interior measurement $u(x)\,(x\in\Omega)$ and the estimates
$||\gamma A||_{L^{\infty}(\Omega)} \le \delta$, $||\rho_k||_{W^{1,\infty}(\Omega)} \le \delta$ $(k=1,2)$ for some constant $\delta>0$.
For $k=1,2$, we denote by $u_k\in H^1(\Omega)$ the solution to
\eqref{eq::fp} with $\rho = \rho_k$ and the constant $\omega^2>0$.
Then as an easy application of the arguments in the previous sections,
we have the following theorem.\\

\begin{theorem}{\rm (Alternate Inverse Problem)}
Let Assumptions 1,2, and 3 be satisfied where $\rho$ is replaced by $\rho_k,\;k=1,2$.  Then, there exist constants $C>0$ and $\alpha\in(0,1)$ depending only on $\Omega$, $d$,
$\sigma$, $g_1$ and the coefficients of $L_1$ such that
\begin{equation}\label{eq::stability_potential}
\Vert\rho_2-\rho_1\Vert_{L^\infty(\Omega_d)}
\le C\,\Vert u_2-u_1\Vert_{W^{1,1}(\Omega)}^\alpha
\end{equation}
for any $g_2$ and $\rho_2$.
Furthermore, if $\Omega$ has an analytic smooth boundary and if
$g_1$ is analytic in $\partial \Omega$ and $\rho_1$, $\gamma A$
are also analytic near $\partial \Omega$, then we have the same estimate
in which $L^\infty(\Omega_d)$ is replaced with $L^\infty(\Omega)$. Note that in this stability estimate \eqref{eq::stability_potential},
we do not have the term $\Vert\rho_2-\rho_1\Vert_{L^\infty(\partial\Omega)}$.
\end{theorem}

\medskip
\begin{proof}
We only point out new considerations that need to be taken in account in applying the arguments in the previous sections.
The key identity we have to use is as follows. For any $\zeta\in H_0^1(\Omega)$,
\begin{equation}\label{eq::identity-potential}
\omega^2\int_\Omega u_1 (\rho_2 - \rho_1) \zeta
=\int_\Omega \nabla_{\gamma A}(u_2-u_1)\cdot\CAnabla\zeta-
\omega^2\int_\Omega\rho_2(u_2-u_1)\zeta.
\end{equation}
Then, by setting $\zeta := \tau_{j,h} \overline{u_1(\rho_2 - \rho_1)}$
in the above key identity,
the proof follows the same arguments as the proof of our Main Theorem. 
\end{proof}

\section{Appendix}
We briefly recall the properties of subanalytic subsets that are needed in our paper.
Reference is made to (\cite{ep}).
Let $X$ and $Y$ be real analytic manifolds.
In what follows, all the manifolds are assumed to be countable at infinity.

\medskip
\begin{definition}{\rm(A subanalytic subset in $X$)} $Z$ is said to be subanalytic at $x\in X$ if there exist
an open neighborhood $U$ of $x$, real analytic compact manifolds $Y_{i,j}\,(i=1,2,\,1\le j\le N)$ and real
analytic maps $\Phi_{i,j}:Y_{i,j}\rightarrow X$ such that
\begin{equation*}
Z\cap U=U\cap\cup_{j=1}^N(\Phi_{1,j}(Y_{1,j})\setminus\Phi_{2,j}(Y_{2,j})).
\end{equation*}
Furthermore, $Z$ is called a subanalytic subset in $X$
if $Z$ is subanalytic at every point $x$ in $X$.
\end{definition}

Let $X$ and $Y$ be real analytic manifolds. The following properties for semi-analytic and subanalytic sets are all found in (\cite{ep}) with their proofs.

\medskip
\begin{enumerate}
\item
Recall that a subset $Z$ in $X$ is said to be semi-analytic if, for any point $x \in X$,
there exists an open neighborhood $V$ of $x$ satisfying
$$
Z \cap V = \underset{i}{\cup}\underset{j}{\cap} \{x \in V\,:\, f_{ij}(x)
\,\,*_{ij} \,\,  0\}
$$
for a finite number of analytic functions $f_{ij}$ on $V$. Here the binary relation
$*_{ij}$ is either $>$ or $=$ for each $i$,$j$.
\item A semi-analytic subset (in particular, an analytic subset) in $X$ is subanalytic in $X$.
\item Let $Z$ be a subset in $X$.  Assume that, for any point $x$
in the closure $\overline{Z}$ of $Z$, there exists
an open neighborhood $V$ of $x$ for which
$Z \cap V$ is subanalytic in $V$. Then $Z$ is subanalytic in $X$.
\item Let $Z$ be a subanalytic subset in $X$. Then its closure, its interior and its
complement in $X$ are again subanalytic in $X$.

\item A finite union and a finite intersection of subanalytic subsets in $X$ are subanalytic in $X$.
\item Let $f: X \to Y$ be a proper analytic map, that is, the inverse image of a compact subset is
again compact. Then, for any subanalytic subset $Z$ in $X$, the image $f(Z)$ is a subanalytic subset
in $Y$.
\end{enumerate}



\medskip
\begin{definition}{\rm (Graph of a subanalytic map)} Let $A$ be a subset in $X$, and let $f: A \to Y$ be a map. We say that
$f$ is a subanalytic map on $A$ if the graph
$$
\Gamma(f):= \{(x,y) \in X \times Y\,:\, x \in A,\, y=f(x)\} \subset X \times Y
$$
is a subanalytic subset in $X \times Y$.
Furthermore, if $Y=\mathbb{R}$,  $f$ is said to be a subanalytic function on $A$.
\end{definition}

\medskip
Note that, if $u$ is a complex valued piecewise analytic function
in an neighborhood of $\overline{\Omega}$ as defined in the body of this paper ,
then $\operatorname{Re}u$, $\operatorname{Im}u$ and $|u(x)|$ are subanalytic functions on $\overline{\Omega}$.
We assume $X = \mathbb{R}^n$ in what follows.
We first recall the following well-known result due to \L{}ojasiewicz
(see Corollary 6.7 in \cite{ep}).\\

\begin{theorem}[\L{}ojasiewicz]\label{Loja}
Let $f$ be a continuous subanalytic function in an open subanalytic subset
$U \subset X = {\mathbb{R}}^n$. Let $Z$ be the zero set of $f$.
For any compact set $K \subset X$, there exist positive constants $C$ and $r$
satisfying
$$
\vert f(x) \vert \ge C \operatorname{dist}(x,\,Z)^r \qquad (x \in U \cap K).
$$
\end{theorem}

\medskip
For the next definition we recall here that a subset $W$ in $X$ is said to be locally closed if
$W$ is a closed subset in an open subset of $X$.
Let $Z$ be a closed subanalytic subset in $X$. \\

\begin{definition}{\rm (Subanalytic Stratification of a Closed Subanalytic subset of X)} We say that a family
$\{Z_\alpha\}_{\alpha \in \Lambda}$ of locally closed subsets
is a subanalytic stratification of $Z$ if the following conditions are satisfied.
\begin{enumerate}
\item $Z$ is a disjoint union of $Z_\alpha$'s.  Each $Z_\alpha$ is called a stratum.
\item $Z_\alpha$ is a connected subanalytic subset in $X$ and
it is analytic smooth at each point in $Z_\alpha$.
\item If $Z_\alpha \cap \overline{Z}_\beta \ne \emptyset$ for $\alpha, \beta \in \Lambda$,
then $Z_\alpha \subset \overline{Z}_\beta$ holds. In particular, we have
$Z_\alpha \subset \partial Z_\beta$ and
$\dim_{\mathbb{R}}Z_\alpha < \dim_{\mathbb{R}}Z_\beta$.
\item The family $\{Z_\alpha\}$ is locally finite in $X$, that is, for any compact
set $K$ in $X$, only a finite number of strata intersect $K$.
\end{enumerate}
\end{definition}

\medskip
For example, let $X=\mathbb{R}^2$
and let us consider a closed triangle $abc$ with its vertexes $a$, $b$ and $c$ as $Z$.
Then $Z$ has a subanalytic stratification
consisting of 7-strata,  the interior of the triangle,
open segments $ab$, $bc$, $ca$ and points $a$, $b$, $c$.  See Figure 1.\\


\begin{figure}[hbt]
\begin{center}
\includegraphics[height=20mm]{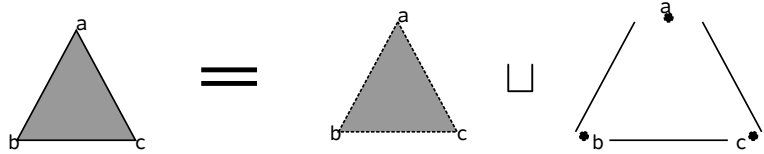}
  \caption{A subanalytic stratification of a closed triangle $Z$.}
\end{center}
\end{figure}

Our interest is in $Z \subset \mathbb{R}^n$ which is a compact subanalytic set 
with $\dim_{\mathbb{R}} Z < n$.
It follows from Theorem A in \cite{kur} that there exists a subanalytic stratification $\{Z_{\alpha}\}_{\alpha \in \Lambda}$ where
each stratum $Z_{\alpha}$ is an L-regular s-cell.  See 6.~Definition in \cite{kur} for the definition of an L-regular s-cell.  Furthermore, since $Z$ is compact, the subanalytic stratum $\{Z_\alpha\}_{\alpha \in \Lambda}$ of $Z$ is locally finite in $X$ implying that the index set $\Lambda$ is finite.

\medskip
The properties of the L-regular s-cell $Z_\alpha$ ($\alpha \in \Lambda$) that we need are that it can be built up from a zero or one dimensional set $B_1$, using orthogonal coordinates $(x_1,...,x_n) \in \mathbb{R}^n$, a positive constant $M$, and where the build up is through ordered pairs $(B_k, \Phi_k)$, $(k=1,...,n-1)$, referred to as data, where $B_k \subset \mathbb{R}^k$ and $\Phi_k$ is a set of functions whose details are given below.  
The stratum, $Z_\alpha$, is thus a kind of cylinder cell built up from a lower dimensional cell to a higher dimensional one; see Figure 2.

\begin{figure}[hbt]
\begin{center}
\includegraphics[height=40mm]{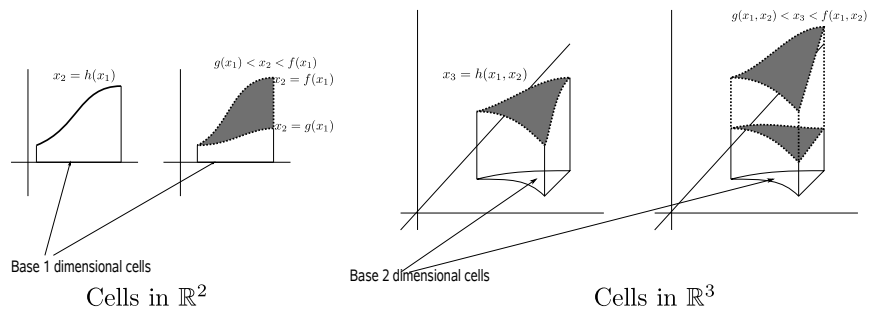}
  \caption{A cylinder cell built up from a base one.}\label{fig:cells}
\end{center}
\end{figure}

\noindent {\it Properties of the L-regular s-cell $Z_\alpha$}

\medskip
\begin{enumerate}
\item[1)] The set $B_1$ is a point or an open interval in $\mathbb{R}^1_{x_1}$.
Each $B_k$ is a locally closed subanalytic subset
in $\mathbb{R}^k_{x_1,\dots,x_k}$ and it is analytic smooth at each point in $B_k$.
\item[2)] The set $\Phi_k$ is a set consisting of one continuous subanalytic function $h_k$ on $\overline{B}_k$
or two continuous subanalytic functions $f_k$ and $g_k$ on $\overline{B}_k$
with $f_k(x') > g_k(x')$ ($x' \in B_k$).
Furthermore, for any $\varphi \in \Phi_k$,
 $\varphi$ is analytic on $B_k$ and has the estimate
\begin{equation}\label{eq:estimate-boundary-function}
\vert d_{B_k} \varphi(x') \vert \le M \qquad (x' \in B_k).
\end{equation}
Here $d_{B_k}\varphi$ denotes the differential 1-form of $\varphi$
on $B_k$ and the cotangent bundle of $B_k$ is
equipped with the metric induced from
the standard one in $\mathbb{R}^k$.
\item[3)] For $k=1,\dots,n-1$,
if $\Phi_k$ consists of one function, then
$$
B_{k+1} = \{(x', x_{k+1}) \in \mathbb{R}^{k+1}\,:\, x' \in B_k,\,
x_{k+1} = h_k(x')\},
$$
where $x' = (x_1,\dots,x_k)$, and otherwise we have
$$
B_{k+1} = \{(x', x_{k+1}) \in \mathbb{R}^{k+1}\,:\, x' \in B_k,\,
g_k(x') < x_{k+1} < f_k(x') \}.
$$
Here we set $B_n := Z_\alpha$.
\end{enumerate}

\medskip
Summing up, as Figure 2
shows,
the L-regular s-cell $Z_\alpha$ is
constructed successively from $B_1$ to $B_n = Z_\alpha$ using functions in $\Phi_1$, \dots, $\Phi_{n-1}$.
Furthermore $Z_\alpha$ itself and each component of $\partial Z_\alpha$
are sufficiently flat due to (\ref{eq:estimate-boundary-function}).

\medskip
Our goal is to present a theorem that gives an open covering, whose measure we can estimate, of the zero level set of a subanalytic function.  Prior to presenting and proving this theorem we establish that we can extend a function in $\Phi_k$ (defined on $\overline B_k \subset \mathbb{R}^k$) to $\mathbb{R}^k$ as a subanalytic Lipschitz continuous function.  We first establish the following lemma.
\\

\begin{lemma} Let $X=\mathbb{R}^n$ and $Z$ be a compact subanalytic subset in $X$.
Let $\{Z_{\alpha}\}_{\alpha \in \Lambda}$ be a subanalytic stratification of $Z$ 
with each $Z_\alpha$ being an L-regular s-cell.  
Further, let $\alpha \in \Lambda$ and $\{(B_k, \Phi_k)\}_{k=1}^{n-1}$ 
be the data for $Z_\alpha$. 
Then $\varphi \in \Phi_k$ is a Lipschitz continuous function on $\overline{B}_k$.
\end{lemma}
\vspace{.1in}

\begin{proof} If we can show the Lipschitz continuity of $\varphi$ on $B_k$,
then the claim of the lemma follows from the continuity of $\varphi$ on $\overline{B}_k$.
Hence it suffices to  prove the claim on $B_k$.
Since $B_k$ itself is an L-regular s-cell in $\mathbb{R}^k$, by 8.~Proposition in \cite{kur},
there exists a positive constant $C$ for which any points $p$ and $q$ in $B_k$ are
joined by a smooth curve $\ell$ in $B_k$ with
\begin{equation}\label{eq:estimate-curve-lenght}
\text{the length of $\ell$} \, \le C|p -q|.
\end{equation}

Let $p$ and $q$ be points in $B_k$, and let $\ell(s)$ ($0 \le s \le 1$) be such a curve in $B_k$.
Then we have
$$
\begin{aligned}
|\varphi(p) - \varphi(q)|
&= \left|\int_0^1 d_{B_k}\varphi\left(\frac{d\ell}{ds}(s)\right) ds\right| \\
&\le M\int_0^1 \left|\frac{d\ell}{ds}(s)\right| ds = M \times (\text{the length of $\ell$}).
\end{aligned}
$$
Here we identified $\frac{d\ell}{ds}(s)$ with a tangent vector of the manifold $B_k$ at $\ell(s)$.
Hence the result follows from (\ref{eq:estimate-curve-lenght}).
\end{proof}

\medskip
Now we construct a family of maps $\rho_k: \mathbb{R}^k \to \overline{B}_k$
($k=1,2,\dots,n-1$) satisfying the following conditions.
\begin{enumerate}
	\item $\rho_k$ is a subanalytic map on $\mathbb{R}^k$ and
		a Lipschitz continuous map on any compact subset in $\mathbb{R}^k$.
	\item $\rho_k(q) = q$ for $q \in \overline{B}_k$.
\end{enumerate}
We construct the family recursively. For $k=1$, we set
$\rho(x) = a$
if $B_1$ consists of one point $a \in \mathbb{R}$, otherwise we define,
for $B_1 = (a,\,b) \subset \mathbb{R}$ ($a < b$),
$$
\rho_1(x) = \left\{
\begin{aligned}
a \qquad &(x < a), \\
x \qquad &(a \le x \le b), \\
b \qquad &(b < x).
\end{aligned}
\right.
$$
Clearly the conditions are satisfied for $\rho_1$. Suppose that $\rho_k$ has been constructed.
We first define $\rho^{(1)}_{k+1}: \mathbb{R}^{k+1} \to \overline{B}_k \times \mathbb{R}_{x_{k+1}}$ by
$$
\rho^{(1)}_{k+1}(x',x_{k+1}) = (\rho_k(x'), x_{k+1}),
$$
which is subanalytic and Lipschitz continuous by the induction hypothesis. 

Now
we define $\rho^{(2)}_{k+1}: \overline{B}_k \times \mathbb{R}_{x_{k+1}}
\to \overline{B}_{k+1}$ in the following way.
If $\Phi_k$ consists of one function $h_k$, we set
$$
\rho^{(2)}_{k+1}(x', x_{k+1}) = (x', h_k(x')).
$$
Otherwise we set
$$
\rho^{(2)}_{k+1}(x', x_{k+1}) = \left\{
\begin{aligned}
(x', g_k(x')) \qquad &(x_{k+1} < g_k(x')), \\
(x', x_{k+1}) \qquad &(g_k(x') \le x_{k+1} \le f_k(x')), \\
(x', f_k(x')) \qquad &(f_k(x') < x_{k+1}).
\end{aligned}
\right.
$$

Since $h_k$ and  $f_k$, $g_k$ are subanalytic
and Lipschitz continuous,
$\rho^{(2)}_{k+1}$ also becomes a subanalytic and Lipschitz continuous map
in both cases.
We set $\rho_{k+1} := \rho^{(2)}_{k+1} \circ \rho^{(1)}_{k+1}$.
Then $\rho_{k+1}$ is a subanalytic and Lipschitz continuous map
as a composition of maps that have
the same properties, and
$\rho_{k+1}(q) = q$ for $q \in \overline{B}_{k+1}$ clearly holds by the construction.
Hence we have obtained the desired family of maps $\rho_k$ $(k=1, \dots, n-1)$.

Let $\varphi \in \Phi_k$. Then $\varphi(\rho_k(x))$ is a subanalytic and
Lipschitz continuous function on
$\mathbb{R}^k$ and its restriction to $\overline{B}_k$ coincides with $\varphi$.
Therefore, in what follows, we assume that all the functions belonging to $\Phi_k$ are
defined in $\mathbb{R}^k$ and they are subanalytic and  Lipschitz continuous there
for any $k=1,2,\dots, n-1$.

It follows from $\operatorname{dim}_{\mathbb{R}} Z_\alpha < n$ that
there exists $1 \le \kappa \le n-1$
such that $\Phi_{\kappa}$ consists of only
one function $h_{\kappa}$.
In fact, otherwise, the $Z_\alpha$ becomes an open subset in $X$
which contradicts $\operatorname{dim}_{\mathbb{R}} Z_\alpha < n$.
Let $\kappa_\alpha$ be the largest one of those $\kappa$'s.
Then we define the subanalytic open subset $U_\alpha(\eta)$ ($\eta > 0$) by
\begin{equation}{\label{eq:def-U_alpha}}
\begin{aligned}
U_\alpha(\eta) =
\{x \in \mathbb{R}^n\,:\,
&h_{\kappa_\alpha}(x_1, \dots, x_{\kappa_\alpha}) - \eta < x_{\kappa_\alpha + 1} <
h_{\kappa_\alpha}(x_1, \dots, x_{\kappa_\alpha}) + \eta, \\
&\qquad \vert x_j \vert < R\,\, (j=1,\dots,\kappa_\alpha,\kappa_{\alpha + 2}, \dots, n) \}.
\end{aligned}
\end{equation}
Clearly $U_\alpha(\eta)$ ($\eta > 0$) is an open subanalytic subset and
it contains $Z_\alpha$.
For the other $\alpha \in \Lambda$, we can construct a subanalytic open neighborhood
$U_\alpha(\eta)$ of $Z_\alpha$ in the same way. 
\\

By setting $U(\eta) = \underset{\alpha \in \Lambda}{\bigcup} U_\alpha(\eta)$ with
$ U_\alpha(\eta)$ defined in the above paragraph, we have the following covering theorem.
\\

\begin{theorem}\label{critical_set}
Let $Z$ be a compact subanalytic subset in $X$ with $\operatorname{dim}_{\mathbb{R}} Z < n$.
Then there exists
a family $U(\eta)$ $($$0 <  \eta \le 1$$)$ of subanalytic open neighborhoods of $Z$
and positive constants $C_1$, $C_2$ for which we have the following.

\begin{enumerate}
\item $\operatorname{vol}(U(\eta)) \le C_1 \eta$ for any $\eta \in (0,\,1]$.
\item
$\operatorname{dist}(p,\, Z) \ge C_2 \eta$ for any point $p \in X \setminus U(\eta)$
and any $\eta \in (0,\,1]$.
\end{enumerate}
\end{theorem}

\medskip
\begin{proof}
We will establish that	$U(\eta)$ has the desired properties described in the statement of the theorem.
Since each $U_\alpha(\eta)$ is subanalytic open and contains $Z_\alpha$,
their union $U(\eta)$ becomes a subanalytic open neighborhood of $Z$.
The first claim 1.~of the theorem is easily seen. In fact, we have
$$
\operatorname{vol}(U_\alpha(\eta)) =
\int_{-R}^{R} \dots \int_{h_{\kappa_\alpha} - \eta}^{h_{\kappa_\alpha} + \eta}\dots
\int_{-R}^{R} dx_1 \dots dx_n =
2(2R)^{n-1}\eta.
$$
Since the number of the strata is finite, the claim follows from this.

We now establish claim 2.~of the theorem. Suppose that the claim were false.
Then there exists a sequence $\{\eta_j\}$ of positive real numbers in $(0,\,1]$
and points $p_j \in X \setminus U(\eta_j)$
satisfying
\begin{equation}\label{eq:contradict-conclusion}
\frac{\operatorname{dist}(p_j,\, Z)}{\eta_j} \to 0 \qquad (j \to \infty).
\end{equation}
Note that, since $\operatorname{dist}(p_j,\, Z) \to 0$ ($j \to \infty$) also holds,
the sequence $\{p_j\}$ is bounded.
Hence, by taking a subsequence, we may assume $\eta_j \to \eta_\infty$ and
$p_j \to p_\infty$ ($j \to \infty$)
for some $\eta_\infty \in [0,1]$ and $p_\infty \in Z$.
Suppose $\eta_\infty > 0$. Then
$p_\infty$ belongs to both $X \setminus U\left(\dfrac{\eta_\infty}{2}\right)$
and $Z$. This
contradicts the fact that $U\left(\dfrac{\eta_\infty}{2}\right)$
is an open neighborhood of the compact
set $Z$. Therefore we assume $\eta_\infty = 0$, i.e., $\eta_j \to 0$ ($j \to \infty$)
in what follows.

Let $q_j$ be a point in $Z$ with $\operatorname{dist}(p_j,\, Z) = |p_j - q_j|$.
By taking a subsequence, we may assume $q_j \in Z_\alpha$ ($j=1,2,\dots$)
for some $\alpha$. Let $\pi_k: \mathbb{R}^n \to \mathbb{R}^k$ ($k=1,2,\dots,n$)
denote the canonical projection defined by
$$
\pi_k(x_1,\dots,x_n) = (x_1, \dots, x_k).
$$
Let $\kappa_\alpha$ be the index determined before equation (\ref{eq:def-U_alpha}).
Then we have
\begin{equation}{\label{eq:xxx_pi_estimate_1}}
| \pi_{\kappa_\alpha}(p_j) - \pi_{\kappa_\alpha}(q_j) | \le
\operatorname{dist}(p_j,\, Z)
\end{equation}
and
\begin{equation}{\label{eq:xxx_pi_estimate_2}}
| \pi_{\kappa_\alpha+1}(p_j) - \pi_{\kappa_\alpha+1}(q_j) | \le
\operatorname{dist}(p_j,\, Z).
\end{equation}
Note that, since $q_j \in Z_\alpha$ and $\Phi_{\kappa_\alpha}$ consists
of only one function $h_{\kappa_\alpha}$, it follows from
the construction of $Z_\alpha$ described
above that
the relation
$$
\pi_{\kappa_\alpha+1}(q_j) =
(\pi_{\kappa_\alpha}(q_j),\,
h_{\kappa_\alpha}(\pi_{\kappa_\alpha}(q_j))) \in \mathbb{R}^{\kappa_\alpha + 1}
$$
holds.

Set
$
\tilde{p}_j = (\pi_{\kappa_\alpha}(p_j),\,
h_{\kappa_\alpha}(\pi_{\kappa_\alpha}(p_j))) \in \mathbb{R}^{\kappa_\alpha + 1}.
$
Then, as $p_j \notin U_\alpha(\eta_j)$
and $\operatorname{dist}(p_j,\, Z_\alpha) \to 0$ ($j \to \infty$),
we have
\begin{equation}{\label{eq:xxx_difference_larger_eta}}
|\tilde{p}_j - \pi_{\kappa_\alpha + 1}(p_j)| \ge \eta_j
\end{equation}
for sufficiently large $j$'s.
Since the function $h_{\kappa_\alpha}$ is Lipschitz continuous, we also have
\begin{equation}{\label{eq:xxx_lipschitz_h}}
\begin{aligned}
|
h_{\kappa_\alpha}(\pi_{\kappa_\alpha}(p_j)) -
h_{\kappa_\alpha}(\pi_{\kappa_\alpha}(q_j)) |
&\le L | \pi_{\kappa_\alpha}(p_j) - \pi_{\kappa_\alpha}(q_j) | \\
&\le L \operatorname{dist}(p_j,\, Z)
\end{aligned}
\end{equation}
for a positive constant $L$. Therefore, by
(\ref{eq:xxx_pi_estimate_1}) and (\ref{eq:xxx_lipschitz_h}),
we obtain
\begin{equation}{\label{eq:xxx_difference_p_and_q}}
|\tilde{p}_j - \pi_{\kappa_\alpha+1}(q_j)|
=
\left|
(\pi_{\kappa_\alpha}(p_j),\, h_{\kappa_\alpha}(\pi_{\kappa_\alpha}(p_j)))
-
(\pi_{\kappa_\alpha}(q_j),\, h_{\kappa_\alpha}(\pi_{\kappa_\alpha}(q_j)))
\right|
\le (1 + L) \operatorname{dist}(p_j,\, Z).
\end{equation}

Summing up, by (\ref{eq:xxx_pi_estimate_2}),
(\ref{eq:xxx_difference_larger_eta}) and
(\ref{eq:xxx_difference_p_and_q}),
we get
$$
\begin{aligned}
\eta_j \le |\tilde{p}_j - \pi_{\kappa_\alpha + 1}(p_j)|
&\le |\tilde{p}_j - \pi_{\kappa_\alpha+1}(q_j)|
+ |\pi_{\kappa_\alpha+1}(q_j) - \pi_{\kappa_\alpha + 1}(p_j)| \\
&\le (2 + L) \operatorname{dist}(p_j,\, Z),
\end{aligned}
$$
from which we have
$$
1 \le (2 + L) \frac{\operatorname{dist}(p_j,\, Z)}{\eta_j}.
$$
This contradicts (\ref{eq:contradict-conclusion}) if $j$ tends to $\infty$, and
hence, the claim 2.~must be true.
The proof has been completed.
\end{proof}

\medskip
\begin{proposition}{\label{prop:area_finite}}
Let $\Omega$ be a relatively compact open subanalytic subset in $\mathbb{R}^n$ and
$f$ a real valued continuous subanalytic function on $\overline{\Omega}$.
Suppose that there exists a subanalytic stratification $\{\Omega_\alpha\}_{\alpha \in \Lambda}$
of $\overline{\Omega}$ such that $f|_{\Omega_\alpha}$ is  analytic in $\Omega_\alpha$
and analytically extends to an open neighborhood of $\overline{\Omega}_\alpha$
for any $\alpha \in \Lambda$ with $\operatorname{dim}_{\mathbb{R}} \Omega_\alpha = n$.
Then there exists a finite subset $E$ of $\mathbb{R}$ and a positive constant
$M_f$ satisfying
$$
H_{n-1}\left(\left\{x \in \overline{\Omega}\,:\, f(x) = t\right\}\right) \le M_f
$$
for any $t \in \mathbb{R}\setminus E$. Furthermore, let $Z$ be a closed subanalytic subset
in $\overline{\Omega}$ with $\operatorname{dim}_{\mathbb{R}} Z \le n-2$. Then
$$
\sup_{t \in \mathbb{R} \setminus E}
H_{n-1}\left(\left\{x \in \overline{\Omega} \cap Z_\epsilon\,
:\, f(x) = t\right\}\right) \to 0 \qquad (\epsilon \to 0^+).
$$
Here $Z_\epsilon := \{x \in \overline{\Omega},\,:\,
\operatorname{dist}(x, Z) \le \epsilon\}$.
\end{proposition}

\

\noindent
{\it{Remark:}}
If $f$ is $C^\infty$, i.e., without subanalyticity, then the claim in the proposition
does not hold even if a subset of measure zero is allowed as $E$.

\

\begin{proof}
For any $t \in \mathbb{R}$, we set $S_t := \{x \in \overline{\Omega}\,:\, f(x) = t\}$.
As we have
$$
\overline{\Omega} =
\underset{\alpha \in \Lambda,\,, \operatorname{dim}_{\mathbb{R}} \Omega_\alpha = n}
{\cup} \overline{\Omega}_\alpha
$$
and $\Lambda$ is a finite set,
it suffices to show the corresponding claim on $\overline{\Omega}_\alpha$
for each $\alpha$ with $\operatorname{dim}_{\mathbb{R}} \Omega_\alpha = n$.
Hence, in what follows, we assume that $f$ is analytic in $\Omega$ and
analytically extendable to an open neighborhood of $\overline{\Omega}$.
If $f$ is a constant function $c$ in $\Omega$, then we take $E=\{c\}$, for which
the claim clearly holds. Therefore we may assume that $f$ is not constant and, as a result,
we have $\operatorname{dim}_{\mathbb{R}} S_t < n$ for any $t$.

Set $\overline{\Omega}_{sing} := \{x \in \overline{\Omega}\,:\, |\nabla f(x)| = 0 \}$
and $\overline{\Omega}_{reg} := \overline{\Omega} \setminus \overline{\Omega}_{sing}$.
Then $f(\overline{\Omega}_{sing})$ is a subanalytic subset in $\mathbb{R}$ as
$f$ is proper on $\overline{\Omega}$ and it is a measure-zero set by Sard's theorem.
Hence $f(\overline{\Omega}_{sing})$ consists of finite points in $\mathbb{R}$
and we take it as $E$.

Let $p_k: \mathbb{R}^n \to \mathbb{R}^{n-1}$ be the canonical projection
by excluding the coordinate $x_k$. We set
$$
\overline{\Omega}_k := \left\{x \in \overline{\Omega}_{reg}\,:\,
|\nabla f(x) - \langle \nabla f(x),\, e_k\rangle e_k| \le n
|\langle \nabla f(x),\, e_k \rangle |\right\},
$$
where $e_k$ is the unit vector with its $k$-th component being $1$.
Note that $\overline{\Omega}_k$ is subanalytic in $\mathbb{R}^n$ and
$\overline{\Omega}_{reg} = \cup_{1 \le k \le n} \Omega_k$
holds. Furthermore we set
$$
S_{t,k} := S_t \cap \overline{\Omega}_k,
$$
which is also subanalytic in $\mathbb{R}^n$.
Then, for any $t \in \mathbb{R} \setminus E$, since $|\nabla f(x)| \ne 0$
on $S_t$, we have $S_t = \cup_{1 \le k \le n} S_{t,k}$ and $S_{t,k}$
is an analytic smooth hypersurface in $\overline{\Omega}_k$.
By these observations it suffices to show
$\underset{t \in \mathbb{R}\setminus E}{\max}\,H_{n-1}(S_{t,k}) < +\infty$.

Define
$$
\begin{aligned}
	\ell_1 &:= \underset{x' \in p_k(\overline{\Omega}_k)}{\max}
	\#\left(\overline{\Omega}_k \cap p_k^{-1}(x')\right), \\
	\ell_2 &:= \underset{x' \in p_k(\overline{\Omega}_k)}{\max}
	\#\left\{x \in \overline{\Omega} \cap p_k^{-1}(x')\,:\,
	\dfrac{\partial f}{\partial x_k} (x) = 0\right\},
\end{aligned}
$$
where $\#A$ denotes the number of the connected components of a set $A$.
Note that these numbers certainly exist because the direct images
$p_{k\,*}\mathbb{R}_{\overline{\Omega}_k}$ and
$p_{k\,*}\mathbb{R}_{\{x \in \overline{\Omega}\,:\, f_{x_k}(x) = 0\}}$
of constructible sheaves
$\mathbb{R}_{\overline{\Omega}_k}$ and
$\mathbb{R}_{\{x \in \overline{\Omega}\,:\, f_{x_k}(x) = 0\}}$
are again constructible sheaves by Proposition 8.4.8 (\cite{ka-sc}) and
$$
\begin{aligned}
\dim_{\mathbb{R}} \left(p_{k\,*}\mathbb{R}_{\overline{\Omega}_k}\right)_{x'}
&= \#\left(\overline{\Omega}_k \cap p_k^{-1}(x')\right),\\
\dim_{\mathbb{R}} \left(p_{k\,*}\mathbb{R}_{\{x \in \overline{\Omega}\,:\, f_{x_k}(x) = 0\}}\right)_{x'}
&= \#\left\{x \in \overline{\Omega} \cap p_k^{-1}(x')\,:\, \dfrac{\partial f}{\partial x_k} (x) = 0\right\},
\end{aligned}
$$
hold for $x' \in p_k(\overline{\Omega}_k)$
(see also Chapter VIII in (\cite{ka-sc})
for the definition of a constructible sheaf).

As $p_k: S_{t,k} \to p_k(S_{t,k})$ is a finite map, that is, $p_k|_{S_{t,k}}$ is a proper map and
$p_k^{-1}(x') \cap S_{t,k}$ consists of finite points for every $x' \in p_k(S_{t,k})$,
there exists a subanalytic stratification $\{O_\beta\}_{\beta \in \Xi}$
of $p_k(S_{t,k})$ such that $S_{t,k} \cap p^{-1}_{k}(O_\beta)$ becomes a
finite covering over $O_\beta$ for each $\beta$.
Note that the stratification consists of a finite number of strata.

Furthermore the number of connected components of $S_{t,k} \cap p^{-1}_{k}(O_\beta)$
is at most $\ell_1 + \ell_2$, which can be proved as follows: As $O_\beta$ is connected, it suffices to show that
the number of points $p^{-1}_k(x') \cap S_{t,k}$ ($x' \in O_\beta$) is at most $\ell_1 + \ell_2$.
Let $L$ be the line $p^{-1}_k(x')$. We first assume that $L \cap \overline{\Omega}_k$ is connected, i.e., $\ell_1 = 1$.
Then there exist mutually distinct points $q_1, \dots, q_m$ in $L$ such that, in each open interval $(q_k,\, q_{k+1})$ of $L$,
$f(x)$ is strictly increasing, strictly decreasing or constant. As $S_{t,k}$ intersects $L$ transversally,
$S_{t,k}$ never intersects an interval where $f(x)$ is constant. Since
the number of intervals in which $f(x)$ is non-constant is at most $\ell_2 + 1$ and
since $S_{t,k}$ intersects
the closure of such an interval at one point if exists,
we conclude that $S_{t,k} \cap L$ consists of at most $\ell_2 +1 = \ell_1 + \ell_2$
points. By applying the same argument to each connected component of
$L \cap \overline{\Omega}_k$, we can prove the claim for the case $\ell_1 > 1$.

For $\beta$ with $\operatorname{dim}_{\mathbb{R}} O_\beta < n -1$,
since $p_k|_{S_{t,k}}$ is a finite covering over $O_\beta$, we have
$$
\operatorname{dim}_{\mathbb{R}}(S_{t,k} \cap p^{-1}_k(O_\beta)) < n-1,
$$
which implies $H_{n-1}(S_{t,k} \cap p^{-1}_k(O_\beta)) = 0$.
On the other hand, for $\beta$ with $\operatorname{dim}_{\mathbb{R}} O_\beta = n -1$,
we have
$$
H_{n-1}(S_{t,k} \cap p^{-1}_k(O_\beta)) \le \sqrt{1 + n^2}(\ell_1 + \ell_2)
H_{n-1}(O_\beta).
$$
Hence we have
\begin{equation}{\label{eq:last-vol-estimate}}
\begin{aligned}
H_{n-1}(S_{t,k})
&\le \sqrt{1 + n^2}(\ell_1 + \ell_2)
\sum_{\beta \in \Xi\,:\, \operatorname{dim}_{\mathbb{R}} O_\beta = n -1}
H_{n-1}(O_\beta) \\
&\le \sqrt{1 + n^2}(\ell_1 + \ell_2) H_{n-1}(p_k(\overline{\Omega})).
\end{aligned}
\end{equation}
This shows the first claim of the proposition.

Finally we show the last claim.  Clearly $\dim_{\mathbb{R}} p_k(Z) < n-1$ and
$$
p_k(Z_\epsilon) \subset \left(p_k(Z)\right)_\epsilon
:= \{y \in \mathbb{R}^{n-1}\,:\, \operatorname{dist}(y, p_k(Z)) \le \epsilon\}
$$
hold.  Hence we have, in $\mathbb{R}^{n-1}$,
$$
0 \le H_{n-1}\left(p_k(Z_\epsilon)\right) \le
H_{n-1} \left(\left(p_k(Z)\right)_\epsilon\right) \to 0, \qquad (\epsilon \to 0^+)
$$
due, for example, to the second Theorem in this Appendix.
Then the last claim of the proposition immediately follows from (\ref{eq:last-vol-estimate}).
\end{proof}


\bibliographystyle{plain}

\end{document}